\documentclass[draft,12pt]{amsart}

\usepackage{amssymb}
\usepackage[curve]{xypic}
\usepackage{amsmath}
\usepackage{amsthm}
\usepackage{amsfonts}
\usepackage{amscd}
\usepackage{pifont}

\newtheorem{Lemma}{Lemma}[section]
\newtheorem{Theorem}[Lemma]{Theorem}
\newtheorem{Corollary}[Lemma]{Corollary}
\newtheorem{Proposition}[Lemma]{Proposition}

\theoremstyle{definition}
\newtheorem{Definition}[Lemma]{Definition}
\newtheorem{Remark}[Lemma]{Remark}

\numberwithin{equation}{section}

\newcommand{\supp}{\operatorname{supp}}
\newcommand{\codim}{\operatorname{codim}}
\newcommand{\Hom}{\operatorname{Hom}}
\newcommand{\Ext}{\operatorname{Ext}}
\newcommand{\rad}{\operatorname{rad}}

\newcommand{\HH}{\operatorname{HH}}
\newcommand{\im}{\operatorname{im}}
\newcommand{\FH}{\mathcal H}
\newcommand{\WH}{\operatorname{WH}}
\newcommand{\OH}{H}

\renewcommand{\AA}{\mathcal A}
\newcommand{\FF}{\mathcal F}
\newcommand{\LL}{\mathcal L}
\newcommand{\QQ}{\mathcal Q}
\newcommand{\CC}{\mathcal C}

\newcommand{\tot}[2]{\buildrel #1 \over #2}

\newcommand{\lcover}{\gtrdot}

\newcommand{\rcover}{\lessdot}
\newcommand{\rcovers}{\rcover}

\allowdisplaybreaks

\begin{document}

\title{The Face Semigroup Algebra of a Hyperplane Arrangement}
\author{Franco V Saliola}
\email{saliola@gmail.com}

\begin{abstract}
This article presents a study of an algebra spanned by the faces of a
hyperplane arrangement. The quiver with relations of the algebra is
computed and the algebra is shown to be a Koszul algebra.
It is shown that the algebra depends only on the intersection lattice of
the hyperplane arrangement.  A complete system of primitive orthogonal
idempotents for the algebra is constructed and other algebraic structure
is determined including: a description of the projective indecomposable
modules; the Cartan invariants; projective resolutions of the simple
modules; the Hochschild homology and cohomology; and the Koszul dual
algebra.  A new cohomology construction on posets is introduced and it is
shown that the face semigroup algebra is isomorphic to the cohomology
algebra when this construction is applied to the intersection lattice of
the hyperplane arrangement.  \end{abstract}

\maketitle

\markboth{\textsc{FRANCO V SALIOLA}}{\textsc{THE FACE SEMIGROUP ALGEBRA}}

\small
\tableofcontents
\normalsize

\section{Introduction}
Let $\AA$ denote a finite collection of linear hyperplanes in $\mathbb
R^d$.  Then $\AA$ dissects $\mathbb R^d$ into open subsets called
\emph{chambers}.  The closures of the chambers are polyhedral cones whose
relatively open faces are called the \emph{faces} of the hyperplane
arrangement $\AA$.  The set $\FF$ of faces of $\AA$ can be endowed with a
semigroup structure. Geometrically, the product $xy$ of faces $x$ and $y$
is the face entered by moving a small positive distance along a straight
line from $x$ towards $y$. The $k$-algebra spanned by the faces of $\AA$
with this multiplication is the \emph{face semigroup algebra} of the
hyperplane arrangement $\AA$. Here $k$ denotes a field.

The face semigroup algebra $k\FF$ has enjoyed recent attention due mainly
to two interesting results. The first result is that a large class of
seemingly unrelated Markov chains can be studied in a unified setting via
the semigroup structure on the faces of a hyperplane arrangement.  The
Markov chains are encoded as random walks on the chambers of a hyperplane
arrangement \cite{BHR1999}. A step in this random walk moves from a chamber
to the product of a face with the chamber according to some probability
distribution on the faces of the arrangement.  This identification
associates the transition matrix of the Markov chain with the matrix of a
linear transformation on the face semigroup algebra of the hyperplane
arrangement.  Questions about the Markov chain can then be answered using
algebraic techniques \cite{Brown2000}. For example, a combinatorial
description of the eigenvalues with multiplicities of the transition matrix
is given and the transition matrix is shown to be diagonalizable. 

The second interesting result concerns the descent algebra of a finite
Coxeter group, a subalgebra of the group algebra of the Coxeter group.  To
any finite Coxeter group is associated a hyperplane arrangement and the
Coxeter group acts on the faces of this arrangement.  This gives an action
of the Coxeter group on the face semigroup algebra of the arrangement.  The
subalgebra of elements invariant under the action of the Coxeter group is
anti-isomorphic to the descent algebra of the Coxeter group
\cite{Bidigare1997, Brown2000}. The descent algebra was introduced in
\cite{Solomon1976} and the proof that it is indeed an algebra is rather
involved.  This approach via hyperplane arrangements provides a new and
somewhat simpler setting for studying the descent algebra. See
\cite{Schocker2005} and \cite{Saliola2006:ReflectionArrangements-arxiv}.

This article presents a study of the algebraic structure of the face
semigroup algebra $k\FF$ of an arbitrary central hyperplane arrangement in
$\mathbb R^d$.  Throughout $k$ will denote a field of arbitrary
characteristic and $\AA$ a finite collection of hyperplanes passing through
the origin in $\mathbb R^d$. The \emph{intersection lattice} of $\AA$ is
the set $\LL$ of intersections of subsets of hyperplanes in $\AA$ ordered
by inclusion.  (Note that some authors order the intersection lattice by
reverse inclusion rather than inclusion.)

The structure of the article is as follows. Sections \ref{section: posets}
and \ref{section: hyperplane arrangements} recall notions from the theory
of posets and hyperplane arrangements, respectively. Section \ref{section:
the face semigroup algebra} defines the face semigroup algebra of a
hyperplane arrangement and describes its irreducible representations.  In
Section \ref{section: primitive idempotents} a complete system of primitive
orthogonal idempotents in $k\FF$ is constructed.  This leads to a
description of the projective indecomposable $k\FF$-modules (Section
\ref{section: projective indecomposable modules}) and a computation of the
Cartan invariants of $k\FF$ (see \ref{cartan invariants}).  The projective
indecomposable modules are used to construct projective resolutions of the
simple $k\FF$-modules in Section \ref{section: projective resolutions of
simple modules}.  The quiver with relations of $k\FF$ is computed in
Section \ref{section: quiver of kF}.  Section \ref{section: the ext-algebra
of kF} proves that $k\FF$ is a Koszul algebra and computes the
$\Ext$-algebra (or Koszul dual) of $k\FF$.  This is used in Section
\ref{section: hochschild cohomology} to compute the Hochschild homology and
cohomology of $k\FF$.  Section \ref{section: connections with poset
cohomology} explores connections with poset cohomology.  A new cohomology
construction is introduced and it is shown that the cohomology algebra,
with its cohomology cup product, is isomorphic to $k\FF$. 
Finally, connections with the Whitney cohomology of the geometric lattice
$\LL^*$ are explored.

\section{Posets}
 \label{section: posets}

This section collects some background from the theory of posets for the
convenience of the reader. An excellent reference is Chapter 3 of
\cite{Stanley1997}.


A \emph{poset} is a finite set $P$ together with a partial order $\leq$.
The \emph{opposite poset} $P^*$ of a poset $P$ is the set $P$ with partial
order defined by $x \leq y$ in $P^*$ iff $x \geq y$ in $P$.  For $x,y \in
P$, write $x \lessdot y$ and say $y$ \emph{covers} $x$ or $x$ is
\emph{covered by} $y$ if $x < y$ and there does not exist $z \in P$ with $x
< z < y$. The \emph{Hasse diagram} of $P$ is the graph with exactly one
vertex for each $x$ in $P$ and exactly one edge between $x$ and $y$ iff $x
\lessdot y$ or $y \lessdot x$. An edge of the Hasse diagram is called a
\emph{cover relation}. 

A \emph{chain} in $P$ is a sequence of elements $x_0 < x_1 < \cdots < x_r$
in $P$.  A chain $x_0 < x_1 < \cdots < x_r$ is \emph{unrefinable} if
$x_{i-1} \rcovers x_{i}$ for all $1 \leq i \leq r$.  The \emph{length} of
the chain $x_0 < x_1 < \cdots < x_r$ is $r$. The \emph{length} or
\emph{rank} of a poset is the length of the longest chain in $P$. For $x
\leq y$ in $P$ the \emph{interval} between $x$ and $y$ is the set $[x,y] =
\{ z \in P \mid x \leq z \leq y\}$. The interval $[x,y]$ is a poset and its
rank is denoted by $\ell([x,y])$.

A (finite) poset $L$ is a \emph{lattice} if every pair of elements $x,y$ in
$L$ has a least upper bound (called \emph{join}) $x \vee y$ and a greatest
lower bound (called \emph{meet}) $x \wedge y$ (with respect to the relation
$\leq$). There exists an element $\hat 0$ called the \emph{bottom} of $L$
satisfying $\hat 0 \leq x$ for all $x \in L$. Similarly, there exists an
element $\hat 1$ in $L$ called the \emph{top} of $L$ satisfying $x \leq
\hat 1$ for all $x \in L$.


The \emph{M\"obius function} $\mu$ of a finite poset $P$ is
defined recursively by the equations
\begin{gather*}
 \mu(x,x) = 1 \quad\text{ and }\quad
 \mu(x,y) = - \sum_{x \leq z < y } \mu(x,z), 
\end{gather*}
for all $x < y$ in $P$. If $x \not<y$, then set $\mu(x,y) = 0$.
The \emph{M\"obius inversion formula} 
\cite[\S3.7]{Stanley1997} states
that
$g(x) = \sum_{y \leq x} f(y)$ 
iff 
$f(x) = \sum_{y \leq x} g(y) \mu(y,x)$,
where
$f,g :P \to \mathbb R$.

\section{Hyperplane Arrangements}
 \label{section: hyperplane arrangements}

This section recalls some background from the theory of hyperplane
arrangements (see \cite{OrlikTerao1992}).

\subsection{Hyperplane Arrangements}

A \emph{hyperplane arrangement} $\AA$ in $\mathbb R^d$ is a finite set of
hyperplanes in $\mathbb R^d$. We restrict our attention to \emph{central}
hyperplane arrangements where all the hyperplanes contain $0 \in \mathbb
R^d$. Each hyperplane $H \in \AA$ determines two open half-spaces of $\mathbb
R^d$ denoted $H^+$ and $H^-$. The choice of which half-space to label $+$
or $-$ is arbitrary, but fixed.

\subsection{The Face Poset}
A \emph{face} of $\AA$ is a nonempty intersection of the form
\begin{gather*}
x = \bigcap_{H\in \AA} H^{\sigma_H},
\end{gather*}
where $\sigma_H \in \{+,-,0\}$ and $H^0 = H$. The sequence $\sigma(x) =
(\sigma_H)_{H\in\AA}$ is the \emph{sign sequence} of $x$. A \emph{chamber}
$c$ is a face such that $\sigma_H(c) \neq 0$ for all $H \in \AA$. 

The \emph{face poset} $\FF$ of $\AA$ is the set of faces of $\AA$ 
partially ordered by 
\begin{gather*}
x \leq y \iff 
\text{for each } H \in \AA \text{ either }
\sigma_H(x) = 0 \text{ or } \sigma_H(x) = \sigma_H(y).
\end{gather*}
Equivalently, $x \leq y \iff x \subset \bar y$.  If $x \leq y$, then we say
$x$ \emph{is a face of} $y$. Note that the chambers are the maximal
elements in this partial order.

\subsection{The Support Map and the Intersection Lattice}
The \emph{support} of a face $x \in \FF$ is the 
the intersection of the hyperplanes in $\AA$ containing $x$,
\begin{gather*}
 \supp(x) = \bigcap_{H \in \AA \atop \sigma_H(x) = 0} H.
\end{gather*}
The set $\LL = \supp(\FF)$ of supports of faces of $\AA$ is a graded lattice
ordered by inclusion, called the \emph{intersection lattice} of $\AA$. (Some
authors order the intersection lattice by \emph{reverse} inclusion, so some
care is needed while reading the literature.) The rank of $X \in \LL$ is
the dimension of the subspace $X \subset \mathbb R^d$ if the intersection
of all the hyperplanes in the arrangement is trivial. For $X,Y\in\LL$
the \emph{meet} $X \wedge Y$ of $X$ and $Y$ is the intersection $X \cap Y$
and the \emph{join} $X \vee Y$ of $X$ and $Y$ is $X+Y$, the smallest subspace
of $\mathbb R^d$ containing $X$ and $Y$. The opposite poset $\LL^*$
of $\LL$ is a geometric lattice. The top element $\hat 1$ of
$\LL$ is the ambient vector space $\mathbb R^d$ and the bottom element $\hat
0$ is the intersection of all hyperplanes in the arrangement $\bigcap_{H \in
\AA} H$. The chambers are the faces of support $\hat 1$. Since $\supp(x) \leq
\supp(y)$ if $x \leq y$, the support map $\supp: \FF \to \LL$ is an
order-preserving poset surjection.

\subsection{Deletion and Restriction}
 \label{restricted arrangements}
Fix $X \in \LL$. The faces $y$ of $\AA$ with $\supp(y) \leq X$ are the
faces of the arrangement $\AA_X = \{ H \cap X \mid X \not\leq H \in \AA \}$.  $\AA_X$
is the \emph{restriction} to $X$ and the face poset of $\AA_X$ is denoted
by $\FF_{\leq X}$.  The intersection lattice $\LL_{\leq X}$ of $\AA_X$ is
the interval $[\hat 0, X]$ of $\LL$.

Given $X \in \LL$ let $\AA^X = \{ H \in \AA \mid X \subset H\}$ denote the
set of hyperplanes in $\AA$ containing $X$. $\AA^X$ is a \emph{deletion} of
$\AA$. If $x \in \FF$ with $\supp(x) = X$, then the face poset $\FF^X$ of
$\AA^X$ is isomorphic to the subposet of $\FF$ of all faces having $x$ as a
face: $\FF^X \cong \{y \in \FF \mid x \leq y\}$. The intersection lattice of
$\AA^X$ is the interval $[X,\hat1] \subset \LL$.

\section{The Face Semigroup Algebra}
 \label{section: the face semigroup algebra}

This section recalls the semigroup structure on the faces of a hyperplane
arrangement and the irreducible representations of the resulting semigroup
algebra. See \cite{Brown2000} for details.

\subsection{The Face Semigroup}
 For $x, y \in \FF$ the product $xy$ is
the face of $\AA$ with sign sequence
\begin{gather*}
\sigma_H(xy) = 
\begin{cases}
\sigma_H(x), & \text{if } \sigma_H(x) \neq 0, \\
\sigma_H(y), & \text{if } \sigma_H(x) = 0.
\end{cases}
\end{gather*}
This product is associative and noncommutative with
identity element the intersection of all the hyperplanes
in the arrangement $1 = \bigcap_{H \in \AA} H$.
Note that the support of the identity element $1$ is $\hat 0$
(and not $\hat 1$).
The support map $\supp: \FF \to \LL$ satisfies 
$
\supp(xy) = \supp(x) \vee \supp(y) \text{ for all } x, y \in \FF.
$
Therefore $\supp$ is a semigroup surjection, where $\LL$ is considered a
semigroup with product given by join $\vee$, as well as an
ordering-preserving poset surjection.

\begin{Remark} There is a nice geometric interpretation of this
product. The face $xy$ is the face that one enters by moving a \emph{small}
positive distance along any straight line from $x$ to $y$. 
\end{Remark}

\begin{Proposition}
 \label{supp properties}
For all $x,y \in \FF$,
\begin{enumerate}
\item\label{LRB identity 1} $x^2 = x$,
\item\label{LRB identity 2} $xyx = xy$,
\item $xy = y$ iff $x \leq y$,
\item\label{xy=x iff XleqY} $xy = x$ iff $\supp(y) \leq \supp(x)$,
\item $\supp(xy) = \supp(x) \vee \supp(y)$,
\end{enumerate}
\end{Proposition}

\begin{Remark} 
Conditions (\ref{LRB identity 1}) and (\ref{LRB identity 2}) 
of the proposition say that 
$\FF$ is a \emph{left regular band}.
\end{Remark}

\subsection{The Face Semigroup Algebra} 

The \emph{face semigroup algebra} of $\AA$ with coefficients in the field
$k$ is the semigroup algebra $k\FF$ of the face semigroup $\FF$ of $\AA$.
Explicitly, it consists of linear combinations of elements of $\FF$ with
multiplication induced by the product of $\FF$. The face semigroup algebra
$k\FF$ is a finite dimensional associative algebra with identity $1 =
\bigcap_{H \in \AA} H$.

Unless explicitly stated otherwise, no assumptions will be made 
on the characteristic of the field $k$.

\subsection{Irreducible Representations}
 \label{irreducible representations}
This section summarizes Section 7.2 of \cite{Brown2000} constructing the
irreducible representations of $k\FF$.

Since $\FF$ and $\LL$ are semigroups, the support map $\supp: \FF \to \LL$
extends linearly to a surjection of algebras $\supp: k\FF \to k\LL$. 
The kernel of this map is nilpotent and the semigroup algebra
$k\LL$ is isomorphic to a product of copies of the field $k$,
one copy for each element of $\LL$. 
This implies that $\ker(\supp)$ is the Jacobson radical of $k\FF$
and that 
the irreducible
representations of $k\FF$ are given by the components of the composition
$k\FF \tot{\supp}\longrightarrow k\LL \tot{\cong}\longrightarrow 
\prod_{X\in\LL} k$. This last map sends $X \in \LL$ to the vector with
$1$ in the $Y$-position if $Y \geq X$ and $0$ otherwise.
The $X$-component of this surjection is
the map $\chi_X: k\FF \to k$ defined on the faces $y \in \FF$ by
\begin{gather*}
 \chi_X(y) = \begin{cases}
  1, & \text{ if } \supp(y) \leq X, \\
  0, & \text{ otherwise}.
 \end{cases}
\end{gather*}
The elements
\begin{equation}\label{idempotents in kL}
E_X = \sum_{Y \geq X} \mu(X,Y) Y,
\end{equation}
one for each $X \in \LL$, correspond to the standard basis vectors of
$\prod_{X\in\LL}k$ under the isomorphism $k\LL \cong \prod_{X \in
\LL} k$ above. They form a basis of $k\LL$ and also form a 
\emph{complete system of
primitive orthogonal idempotents} (see Section \ref{section: primitive
idempotents}).

\section{Primitive Idempotents}
 \label{section: primitive idempotents}

Let $A$ be a $k$-algebra. An element $e \in A$ is \emph{idempotent} if $e^2
= e$.  It is a \emph{primitive idempotent} if $e$ is idempotent and we
cannot write $e = e_1 + e_2$ where $e_1$ and $e_2$ are nonzero idempotents in
$A$ with $e_1e_2 = 0 = e_2e_1$.  Equivalently, $e$ is primitive iff $A e$
is an indecomposable $A$-module. A set of elements $\{e_i\}_{i\in I} \subset A$ is
a \emph{complete system of primitive orthogonal idempotents} for $A$ if $e_i$
is a primitive idempotent for every $i$, if $e_ie_j = 0$ for $i \neq j$ and
if $\sum_i e_i = 1$. If $\{e_i\}_{i\in I}$ is a complete system of primitive
orthogonal idempotents for $A$, then $A \cong \bigoplus_{i\in I} Ae_i$ as left
$A$-modules and $A \cong \bigoplus_{i,j\in I} e_iAe_j$ as $k$-vector spaces.

\subsection{Complete System of Primitive Orthogonal Idempotents}
 \label{primitive idempotents}
For each $X \in \LL$, fix an $x \in \FF$ with $\supp(x) =
X$ and define elements in $k\FF$ recursively by the formula,
\begin{gather}
 \label{equation: idempotents}
 e_X = x - \sum_{Y > X} x e_Y.
\end{gather}
Note that $e_{\hat1}$ is an arbitrarily chosen chamber.

\begin{Lemma} \label{idempotent lemma}
 Let $w \in \FF$ and $X \in \LL$. If $\supp(w) \not\leq X$, then $w e_X = 0$.
\end{Lemma}
\begin{proof}
We proceed by induction on $X$. This is vacuously true if $X = \hat 1$.
Suppose the result holds for all $Y \in \LL$ with $Y >
X$. Suppose $w \in \FF$ and $W = \supp(w) \not\leq X$. Using
the definition of $e_X$ and the identity $wxw = wx$
(Proposition \ref{supp properties} (\ref{LRB identity 2})),
\begin{align*}
w e_X & = wx - \sum_{Y > X} wx e_Y = wx - \sum_{Y > X} wx(we_Y).
\end{align*}
By induction, $we_Y = 0$ if $W \not\leq Y$. Therefore, the
summation runs over $Y$ with $W \leq Y$. 
But $Y > X$ and $Y \geq W$ iff $Y \geq W \vee X$, so the summation runs
over $Y$ with $Y \geq W \vee X$. 
\begin{align*}
we_X = wx - \sum_{Y > X} wx(we_Y) = wx - \sum_{Y \geq X \vee W} wx e_Y.
\end{align*}
Now let $z$ be the element of support $X \vee W$ chosen in defining $e_{X
\vee W}$. So $e_{X\vee W} = z - \sum_{Y>X\vee W} ze_Y$. 
Note that $ze_{X\vee W} = e_{X\vee W}$ since $z = z^2$. Therefore, $z =
\sum_{Y\geq X\vee W} ze_Y$. Since $\supp(wx) = W \vee X = \supp(z)$, it
follows from 
Proposition \ref{supp properties} (\ref{xy=x iff XleqY})
that
$wx = wxz$.
 Combining the last two statements,
\begin{align*}
we_X = wx - \sum_{Y \geq X \vee W} wx e_Y 
& = wx\left(z - \sum_{Y \geq X \vee W} z e_Y\right) = 0.
\qedhere
\end{align*}
\end{proof}

\begin{Theorem}
 \label{complete system of primitive orthogonal idempotents}
The elements $\{e_X\}_{X \in \LL}$ form a complete
system of primitive orthogonal idempotents in $k\FF$.
\end{Theorem}
\begin{proof}
\emph{Complete}. $1 = \bigcap_{H\in\AA} H$ is the only element of support
$\hat 0$. Hence, $e_{\hat 0} = 1 - \sum_{Y > \hat 0} e_Y$. 
Equivalently, $1 = \sum_{Y \in \LL} e_Y$.

\emph{Idempotent.} 
Since $e_Y$ is a linear combination of elements of support at least $Y$, $e_Y
z = e_Y$ for any $z$ with $\supp(z) \leq Y$
(Proposition \ref{supp properties} (\ref{xy=x iff XleqY})). 
Using the definition of $e_X$,
the facts $e_X = xe_X$ and $e_Y = e_Yy$, and Lemma \ref{idempotent lemma},
\begin{align*}
e_X^2 = \left(x - \sum_{Y > X} xe_Y\right)e_X =
xe_X - \sum_{Y > X} xe_Y(ye_X) = xe_X = e_X.
\end{align*}

\emph{Orthogonal.} We show that for every $X \in \LL$,
$e_Xe_Y = 0$ for $Y \neq X$. If $X = \hat 1$, then $e_X e_Y = e_X x e_Y = 0$
for every $Y \neq X$ by Lemma \ref{idempotent lemma} since $X = \hat 1$
implies $X \not\leq Y$.  Now suppose the result holds for $Z > X$. That is,
$e_Ze_Y = 0$ for all $Y \neq Z$. If $X \not\leq Y$, then $e_Xe_Y = 0$ by
Lemma \ref{idempotent lemma}. If $X < Y$, then $e_X e_Y = xe_Y - \sum_{Z > X}
x (e_Ze_Y) = xe_Y - xe_Y^2 = 0.$

\emph{Primitive.} We'll show that $e_X$ lifts $E_X = \sum_{Y\geq
X} \mu(X,Y)Y$ (see equation (\ref{idempotents in kL})) for all $X \in \LL$,
a primitive idempotent in
$k\LL$. If $X = \hat 1$, then $\supp(e_{\hat 1}) = \hat 1 =
E_{\hat 1}$. Suppose the result holds for $Y > X$. Then
$\supp(e_X) = \supp(x - \sum_{Y>X}xe_Y) = X -
\sum_{Y>X}(X\vee E_Y)$. Since $E_Y$ is a linear combination
of elements $Z \geq Y$, it follows that $X \vee E_Y = E_Y$ if $Y>X$.
Therefore, $\supp(e_X) = X - \sum_{Y>X}E_Y$. The M\"obius
inversion formula applied to $E_X = \sum_{Y\geq X} \mu(X,Y)Y$
gives $X = \sum_{Y \geq X} E_X$. Hence,
$\supp(e_X) = X - \sum_{Y>X}E_Y = E_Y$.

To see that this is sufficient, suppose $E$ is a primitive idempotent in
$k\LL$ and that $e$ is an idempotent lifting $E$. Suppose $e = e_1 + e_2$
with $e_i$ orthogonal idempotents. Then $E = \supp(e) = \supp(e_1) +
\supp(e_2)$. Since $E$ is primitive and $\supp(e_1)$ and $\supp(e_2)$ are
orthogonal idempotents, $\supp(e_1) = 0$ or $\supp(e_2) = 0$. Say
$\supp(e_1) = 0$. Then $e_1$ is in the kernel of $\supp$. This kernel is
nilpotent so $e_1^n = 0$ for some $n \geq 0$. Hence $e_1 = e_1^n = 0$. 
Therefore, $e$ is a primitive idempotent.
\end{proof}

\begin{Remark}
 \label{Remark: different idempotents}
 We can replace $x \in \FF$ in 
 Equation (\ref{equation: idempotents}) with any linear
 combination  $\tilde x = \sum_{\supp(x)=X} \lambda_x x$ 
 of elements of support $X$ whose
 coefficients $\lambda_x$ sum to $1$.
The proofs still hold since the element $\tilde x$ 
is idempotent 
and satisfies
$\supp(\tilde x) = X$ and $\tilde x y = \tilde x$
if $\supp(y) \leq X$. 
Unless explicitly stated we will use the idempotents constructed
above.
\end{Remark}

\subsection{A Basis of Primitive Idempotents}

\begin{Proposition} \label{basis of idempotents}
The set $\{ xe_{\supp(x)} \mid x \in \FF\}$ is a basis of $k\FF$
of primitive idempotents.
\end{Proposition}
\begin{proof}
Let $y \in \FF$. Then by Corollary \ref{complete system of primitive
orthogonal idempotents} and Lemma \ref{idempotent lemma},
$$
y = y1 = y \sum_X e_X = \sum_{X \geq \supp(y)} ye_X = \sum_{X \geq \supp(y)}
(yx)e_X.
$$
Since $\supp(yx) = \supp(y) \vee \supp(x) = X$, the face $y$ is a linear
combination of the elements of the form $xe_{\supp(x)}$. So these elements
span $k\FF$. Since the number of these elements is the cardinality of $\FF$,
which is the dimension of $k\FF$, the set forms a basis of $k\FF$. The
elements are idempotent since $(x e_X)^2 = (x e_X)(x e_X) = x e_X^2 = x e_X$
(since $xyx=xy$ for all $x,y\in\FF$). Since $xe_X$ also lifts the
primitive idempotent $E_X = \sum_{Y\geq X} \mu(X,Y)Y \in k\LL$,
it is also a primitive idempotent (see the end of the proof of
Corollary \ref{complete system of primitive orthogonal idempotents}).
\end{proof}

\section{Projective Indecomposable Modules}
 \label{section: projective indecomposable modules}

This section describes the projective indecomposable $k\FF$-modules
and computes the Cartan invariants of $k\FF$.

\subsection{Projective Indecomposable Modules}

For $X \in \LL$, let $\FF_X \subset \FF$ denote the
set of faces of support $X$. For $y \in \FF$ and $x \in \FF_X$
let
\begin{gather*}
y\cdot x = \begin{cases}
yx, & \supp(y) \leq \supp(x), \\
0,& \supp(y) \not\leq \supp(x).
\end{cases}
\end{gather*}
Then $k\FF_X$ is a $k\FF$-module.

\begin{Lemma}
 \label{basis for indecomposable modules}
Let $X \in \LL$. Then $\{y e_X \mid \supp(y) = X\}$ is a basis
for $k\FF e_X$.
\end{Lemma}
\begin{proof}
Suppose $\sum_{w \in \FF} \lambda_w we_X \in k\FF e_X$. 
If $\supp(w) \not\leq X$, then $we_X = 0$. So suppose $\supp(w) 
\leq X$. Then $\supp(wx) = \supp(w) \vee X = X$.
Therefore,
$$
\sum_{w \in \FF} \lambda_w we_X = \sum_{w \in \FF} 
\lambda_w (wx) e_X \in \operatorname{span}_k\{y e_X \mid \supp(y) = X\},
$$
where $x$ is the element chosen in the construction of $e_X$ (recall
that $e_X = xe_X$ since $x^2 = x$).
So the elements span $k\FF e_X$. These elements are linearly independent
being a subset of a basis of $k\FF$ (Proposition \ref{basis of idempotents}).
\end{proof}

\begin{Proposition}
 \label{projective indecomposable modules}
The $k\FF$-modules $k\FF_X$ are all the
projective indecomposable $k\FF$-modules. The radical
of $k\FF_X$ is $\operatorname{span}_k\{y - y' \mid y,y' \in \FF_X\}$.
\end{Proposition}
\begin{proof}
Define a map $\phi: k\FF_X \to k\FF e_X$ by $w \mapsto
we_X$.  Then $\phi$ is surjective since $\phi(y) = ye_X$
for $y \in \FF_X$ and since $\{ye_X \mid \supp(y) = X\}$
is basis for 
$k\FF e_X$ (Lemma \ref{basis for indecomposable modules}).
Since $\dim k\FF_X = \#\FF_X = \dim k\FF e_X$,
the map $\phi$ is an isomorphism of $k$-vector spaces.

\emph{$\phi$ is a $k\FF$-module map}. Let $y \in \FF$
and let $x \in \FF_X$. If $\supp(y) \leq X$, then $\phi(y \cdot
x) = \phi(yx) = yxe_X = y \phi(x)$. If $\supp(y) \not\leq X$, then $y
\cdot x = 0$. Hence, $\phi(x \cdot y) = 0$. Also, since $\supp(y)
\not\leq X$, it follows that $y e_X = 0$. Therefore, $y\phi(x) = y xe_X = yx(y
e_X) = yx 0 = 0.$ So $\phi(y\cdot x) = y \phi(x)$. Hence $\phi$ is an
isomorphism of $k\FF$-modules. Since $k\FF e_X$ are all the projective
indecomposable $k\FF$-modules, so are the $k\FF_X$.
\end{proof}

\subsection{Cartan Invariants}

Let $\{e_X\}_{X \in I}$ be a complete system of primitive
orthogonal idempotents for a finite dimensional $k$-algebra $A$.
The \emph{Cartan invariants} of $A$ are defined to be the numbers
$$
 c_{X,Y} = \dim \Hom_{A}(A e_X, A e_Y),
$$
where $X,Y \in I$. The invariant $c_{X,Y}$ is the multiplicity of the
simple module $S_X = (A / \text{rad}A)e_X$ as a composition factor of the
left $A$-module $A e_Y$. The \emph{Cartan matrix} of $A$ is the matrix
$[c_{X,Y}]$.

The following is Theorem $1.7.3$ of \cite{Benson1998:I}.
\begin{Theorem}[Idempotent Refinement Theorem]
 \label{idempotent refinement theorem}
Let $N$ by a nilpotent ideal in a ring $R$ and let $e$ be an idempotent
in $R/N$. Then any two idempotents in $R$ lifting $e$ are conjugate in $R$.
\end{Theorem}

\begin{Proposition}
 \label{cartan invariants}
 For $X,Y \in \LL$, 
 $$
 \dim_k\Hom_{k\FF}(k\FF e_X, k\FF e_Y) = |\mu(X,Y)|,
 $$
 where $\mu$ is the M\"obius function of $\LL$. Therefore the Cartan
 invariants of $k\FF$ are $c_{X,Y} = |\mu(X,Y)|$ and the Cartan matrix is
 triangular of determinant $1$.
\end{Proposition}
\begin{proof}
Since $\Hom_{k\FF}(k\FF e_X, k\FF e_Y) \cong e_X k\FF e_Y$, it follows
that $c_{X,Y} = \dim
e_X k\FF e_Y$. We will use Zaslavsky's Theorem \cite{Zaslavsky1975}:
The number of chambers in a hyperplane arrangement is $\sum_{X \in \LL}
|\mu(X,\mathbb R^d)|$.

For each $W \in \LL$, let $w$ denote an element of support $W$. If $W \geq
X$, then $\supp(xw) = W$, so replace $w$ with $xw$ and construct idempotents
$e_W$ as in section \ref{primitive idempotents}. (By the \emph{idempotent
refinement theorem} above, it does not matter
which lifts of the idempotents in $k\LL$ we use to compute the Cartan
invariants: $e_X k\FF e_Y \cong \tilde e_X k\FF \tilde e_Y$ if
$e_X$ and $\tilde e_X$ are conjugate and if $e_Y$ and $\tilde e_Y$
are conjugate.) Then for each $W \geq X$ we have $x e_W = e_W$, so $x = x
\sum_W e_W = x \sum_{W \geq X} e_W = \sum_{W \geq X} e_W$. This gives the
equality 
\begin{gather}
 \label{deleted face semigroup algebra identity}
 k(x\FF) = xk\FF = \sum_{W \geq X} e_W k\FF.
\end{gather}
Note that $x\FF$ is the face poset of the hyperplane arrangement $\AA^X = \{ H \in \AA
\mid X \subset H \}$ and that the faces of 
support $Y$ in $\AA^X$ are the chambers in the restricted arrangement
$(\AA^X)_Y$ (see Section \ref{restricted arrangements}).
Zaslavsky's Theorem
applied to $(\AA^X)_Y$ gives the number of faces of support $Y$ in $\AA^X$ is
$\sum_{W \in [X,Y]} |\mu(W,Y)|$
since the intersection lattice of $(\AA^X)_Y$ is the interval $[X,Y]$ in $\LL$.
But the number of faces of support $Y$ in
$(\AA^X)_Y$ is the cardinality of the set $x\FF_Y$, which is the dimension of
$k(x\FF_Y) \cong x k\FF_Y \cong x k\FF e_Y \cong \bigoplus_{X \leq W \leq Y}
e_W k\FF e_Y$ by (\ref{deleted face semigroup algebra identity}) and
Lemma \ref{idempotent lemma}.
Therefore for each $X, Y \in \LL$,
$$
 \sum_{X \leq W \leq Y} \dim e_W k\FF e_Y = \sum_{X \leq W \leq Y} |\mu(W,Y)|.
$$
The result now follows by induction. If $X = Y$,
then $\dim e_X k\FF e_X = |\mu(X,X)|$. Suppose the result
holds for all $W$ with $X < W \leq Y$. Then
\begin{align*}
 \dim e_X k\FF e_Y
 &= \sum_{X \leq W \leq Y} |\mu(W,Y)| - \sum_{X < W \leq Y} \dim e_W k\FF e_Y \\
 &= \sum_{X \leq W \leq Y} |\mu(W,Y)| - \sum_{X < W \leq Y} |\mu(W,Y)| \\
 &= |\mu(X,Y)|. 
\qedhere
\end{align*}
\end{proof}

\section{Projective Resolutions of the Simple Modules}
 \label{section: projective resolutions of simple modules}

\subsection{A Projective Resolution of the Simple Module 
  Corresponding to $\hat1$}

In Section 5C of \cite{BrownDiaconis1998} an exact sequence of $k\FF$-modules is
constructed to compute the multiplicities of the eigenvalues of random walks
on the chambers of a hyperplane arrangement. This construction
in combination with the above description of the projective indecomposable
$k\FF$-modules
yields a projective resolution of the simple $k\FF$-modules.

Let $\FF_p \subset \FF$ denote the set of faces of 
codimension $p$. For $x \in \FF$ and $y \in \FF_p$, let
\begin{gather*}
x\cdot y = \begin{cases}
xy, & \supp(x) \leq \supp(y), \\
0,& \supp(x) \not\leq \supp(y).
\end{cases}
\end{gather*}

Fix an orientation $\epsilon_X$ for every subspace $X \in \LL$. If $x$ is a
codimension one face of $y$, then pick a positively oriented basis $\{e_1,
\ldots, e_i\}$ of $X = \supp(x)$ and a vector $v$ in $y$ and put
$$
[x:y] = \epsilon_Y(e_1, \cdots, e_i, v),
$$
where $Y = \supp(y)$. Since $X$ is a codimension one subspace of $Y$, the
mapping $v
\mapsto \epsilon_Y(e_1, \cdots, e_i, v)$ is constant on the open
halfspaces of $Y$ determined by $X$. This implies the identity,
\begin{gather} \label{incidence number identity}
 [x:y] = [\tilde x: \tilde x y], \text{ if } \supp(\tilde x) = \supp(x).
\end{gather}
\begin{Lemma}[\cite{BrownDiaconis1998}, \S5 Lemma 2] 
 \label{incidence number lemma}
Let $x, y \in \FF$ with $x$ of codimension two in $y$. Then
there are exactly two faces $w$ and $z$ in the open interval
$(x,y)$ and we have
$$
[x:w][w:y] = - [x:z][z:y].
$$
\end{Lemma}

\begin{Proposition}[\cite{BrownDiaconis1998}, \S5 Lemma 4]
 \label{Brown-Diaconis exact sequence}
The following is an exact sequence of $k\FF$-modules.
$$
\xymatrix{
 \cdots \ar[r] &
 k\FF_p \ar[r]^{\partial_p} &
 \cdots \ar[r] &
 k\FF_1 \ar[r]^{\partial_1} &
 k\FF_0 \ar[r]^{\partial_0} &
 k \ar[r] &
 0, 
}
$$
where the action of $k\FF$ on $k$ is given by 
$w \cdot \lambda = \lambda$ for all $w \in \FF$ and $\lambda \in k$.
The differential $\partial_i$ is given by
$\partial_0(c) = 1$ for all $c \in \FF_0$ and for $x \in \FF_p$,
$$
 \partial_p(x) = \sum_{y \lcover x} [x: y] y.
$$
\end{Proposition}

\begin{proof}[Sketch of the proof]
It is easy to check that the complex consists of $k\FF$-modules
and that $\partial_i$ is a $k\FF$-module map. It remains to explain
why the complex is exact.
Suppose that the intersection of all the hyperplanes in a point,
otherwise quotient out by that subspace.
Intersecting the hyperplane arrangement with a sphere centered
at the origin induces a regular cell decomposition $\Sigma$ of the
$(d-1)$-sphere whose cells correspond to the faces $x \neq 1$ of $\AA$. The
dual of $\Sigma$ is the boundary of a polytope (a zonotope, actually) $Z$.
Therefore, the poset of nonempty faces of $Z$ is anti-isomorphic to the
face poset $\FF$ of $\AA$. Since $Z$ is contractible any augmented cellular
chain complex will be an exact sequence of $k$-vector spaces.
The above complex is precisely the augmented cellular chain complex
with incidence numbers given by $[x:y]$. 
(See \cite{CookeFinney1967}.)
Therefore, it is exact.
\end{proof}

Note that $k\FF_p \cong \bigoplus_{\codim(X) = p} k \FF_X$ as
$k\FF$-modules and that $k\FF_X$ is projective by Proposition
\ref{projective indecomposable modules}, where $\codim(X)$ is the 
codimension of the subspace $X$. So the $k\FF$-modules $k\FF_p$ are
projective. Also note that in order for $\partial_0$ to be a $k\FF$-module
morphism, the action of $k\FF$ on $k$ must be given by $\chi_{\hat 1}$.
That is, $k$ is the simple module afforded by the irreducible
representation $\chi_{\hat 1}$. This proves the following result.

\begin{Corollary}
The exact sequence
$$
\xymatrix{
 \cdots \ar[r] &
 k\FF_p \ar[r]^{\partial_p} &
 \cdots \ar[r] &
 k\FF_1 \ar[r]^{\partial_1} &
 k\FF_0 \ar[r]^{\partial_0} &
 k \ar[r] &
 0
}
$$
is a projective resolution of the simple $k\FF$-module afforded
by the irreducible representation $\chi_{\hat1}: k\FF \to k$.
\end{Corollary}

\subsection{Projective Resolutions of the Simple Modules}

Recall that the simple $k\FF$-modules 
are indexed by $X \in \LL$, afforded by the representations
$\chi_X: k\FF \to k$,
\begin{gather*}
 \chi_X(y) = \begin{cases}
  1, & \text{ if } \supp(y) \leq X, \\
  0, & \text{ otherwise}.
 \end{cases}
\end{gather*}
Also recall that $\FF_{\leq X}$ denotes the face semigroup of $\AA_X$, 
consisting
of the set of faces in $\FF$ of support contained in $X$ (Section
\ref{restricted arrangements}). Let $(\FF_{\leq X})_p$ denote the set of
faces in $\AA_X$ of codimension $p$ in $X$. Applying the previous result to
the hyperplane arrangement $\AA_X$ gives a projective resolution
$$
 \cdots \longrightarrow
 k(\FF_{\leq X})_p \tot{\partial}\longrightarrow
 \cdots \tot{\partial}\longrightarrow
 k(\FF_{\leq X})_1 \tot{\partial}\longrightarrow
 k \FF_X \longrightarrow 
 k_X \longrightarrow
 0
$$
of the simple $k\FF_{\leq X}$-module $k_X$ with action given by
$w \cdot \lambda = \lambda$ for
all $w \in \FF_{\leq X}$ and $\lambda \in k$. The algebra surjection $k\FF
\to k\FF_{\leq X}$ given by $w \mapsto \chi_X(w)w$ for $w \in \FF$ puts a
$k\FF$-module structure on each $k(\FF_{\leq X})_p$ and on $k$. The
$k\FF$-module structure on $k$ is precisely that given by $\chi_X:k\FF \to
k$. Each $k(\FF_{\leq X})_p$ is a projective $k\FF$-module since the
$k\FF$-module structure on $k(\FF_{\leq X})_p$ decomposes as
$$
k(\FF_{\leq X})_p \cong \bigoplus_{{Y \leq X,} \atop {\codim_X(Y)=p}} 
k\FF_Y,
$$
where $\codim_X(Y)$ denotes the codimension
of $Y$ in $X$. This establishes the following.

\begin{Proposition} \label{projective resolution}
Let $X \in \LL$. Then
$$
 \cdots \longrightarrow
\left(
 \bigoplus_{Y \in \LL \atop \codim_X(Y) = p} k \FF_Y
\right) \tot{\partial}\longrightarrow
 \cdots \tot{\partial}\longrightarrow
 k \FF_X \longrightarrow 
 k_X \longrightarrow
 0
$$
is a projective
resolution of the simple $k\FF$-module $k_X$ afforded
by $\chi_X:k\FF \to k$,
where $\partial(w) = \sum_{y \lcover w} [w: y] \chi_X(y) y$
and $\codim_X(Y)$ denotes the codimension of $Y$ in $X$.
\end{Proposition}

\section{The Quiver of the Face Semigroup Algebra}
 \label{section: quiver of kF}

\subsection{The Quiver of a Split Basic Algebra}
A finite dimensional
$k$-algebra $A$ is a \emph{(split) basic algebra} if every simple
module of $A$ has dimension one. The \emph{$\Ext$-quiver} 
or just \emph{quiver} $Q$ of a split basic algebra $A$ is a directed graph with
one vertex for each isomorphism class of simple modules of $A$.
The number of arrows $x \to y$ is $\dim \Ext^1_A(S_x, S_y)$, where $S_x$ and
$S_y$ are simple modules corresponding to the vertices $x$ and
$y$. 

A \emph{path} $p$ in $Q$ is a sequence of arrows $x_0 \to x_1 \to \cdots \to
x_r$. The path \emph{starts} at $s(p) = x_0$ and \emph{terminates} at $t(p) =
x_r$. The length of $p$ is $r$.
Two paths $p$ and $q$ are \emph{parallel} if they start and terminate at the
same vertices: $s(p) = s(q)$ and $t(p) = t(q)$. The \emph{path algebra} $kQ$
of a quiver $Q$ is the $k$-vector space spanned by the paths in $Q$ with the
product of two paths defined by \emph{path composition}: if 
$p = x_0 \to x_1 \to \cdots \to x_r$ and
$q = y_0 \to y_1 \to \cdots \to y_s$, then
$$
p \cdot q =
\begin{cases}
y_0 \to \cdots \to y_s \to x_1 \to \cdots \to x_r, & \text{if } x_0 = s(p)
= t(q) = y_s, \\
0, & \text{otherwise}.
\end{cases}
$$

Let $P \subset kQ$ be the ideal of $kQ$ generated by the arrows of $Q$.
An ideal $I \subset kQ$ is \emph{admissible} if $P^r \subset I \subset
P^2$, for some $r \geq 2$.
\begin{Proposition}[\cite{ARS1995}, \S III.1 Thereom 1.9]
 \label{main result about quivers}
Let $A$ be a finite dimensional split basic $k$-algebra with quiver
$Q$. Then $A \cong kQ/I$ where $I$ is an admissible ideal of $kQ$.
\end{Proposition}

Let $I$ be an admissible ideal of $kQ$. An element of $I$ is a
\emph{relation} from $x$ to $y$ if it is a $k$-linear combination of paths
in $Q$ beginning at a vertex $x$ and ending at a vertex $y$.  Note that any
element $\rho \in I$ can be written as a linear combination of relations
since $x \rho y$ is a relation for any pair of vertices $x, y \in Q$. The
following result combines Corollary $1.1$ and Proposition $1.2$ of
\cite{Bongartz1983}. 


%

\begin{Proposition}\label{minimal relations}
Let $Q$ be a quiver with no oriented cycles and let $I$ be an admissible
ideal. Suppose that $R$ is a minimal set of relations generating $I$ as a
two-sided ideal of $kQ$.  Then the number of relations from $x$ to
$y$ in $R$ is the dimension of the $k$-vector space $\Ext^2_{kQ/I}(S_x, S_y)$.
\end{Proposition}

\subsection{The Quiver of the Face Semigroup Algebra}
Since every simple $k\FF$-module is of dimension one, $k\FF$ is a split
basic algebra. This section computes the quiver $\QQ$ of $k\FF$
and the next section describes an ideal $I$ such that $k\QQ/I \cong
k\FF$.

\begin{Lemma} 
 \label{ext-spaces}
For $X, Y \in \LL$ and $p \geq 0$,
$$
\Ext^p_{k\FF}(k_X, k_Y) \cong
\begin{cases}
k, & \text{if } Y \leq X \text{ and } \dim(X) - \dim(Y) = p, \\
0, & \text{otherwise}.
\end{cases}
$$
\end{Lemma}
\begin{proof}
Let $\codim_X(W)$ denote the codimension of $W$ in $X$
and let $C_p$ denote $\bigoplus_{\codim_X(W)=p} k\FF_W$. Applying the functor
$\Hom(-, k_Y)$ to the projective resolution of $k_X$ in Proposition
\ref{projective resolution}, gives the cocomplex
$$
\cdots \tot{\partial_p^*}\longrightarrow
\Hom_{k\FF}\left(C_p, k_Y\right)
\tot{\partial_{p+1}^*}\longrightarrow 
\Hom_{k\FF}\left(C_{p+1}, k_Y\right)
\tot{\partial_{p+2}^*}\longrightarrow 
\cdots.
$$
Now
$
 \Hom_{k\FF}\left(C_p, k_Y\right)
 \cong
 \bigoplus_{\codim_X(W)=p} \Hom_{k\FF}\left(k\FF_W, k_Y\right)
$
\ and
$$
\Hom_{k\FF}\left(k\FF_W, k_Y\right) \cong
\Hom_{k\FF}\left(k\FF e_W, k_Y\right) \cong
e_W \cdot k_Y = 
\begin{cases}
k, & \text{if } W = Y, \\
0, & \text{otherwise},
\end{cases}
$$
where we used the fact that $\chi_Y(e_W) = 0$ if $W \neq Y$ and $1$
otherwise. (If $W \neq Y$, then $\chi_Y(e_Y) = 1$ implies $\chi_Y(e_W) =
\chi_Y(e_W)\chi_Y(e_Y) = \chi_Y(e_We_Y) = 0$.)
%
%
Since $\Hom_{k\FF}(k\FF_W,k_Y)$ vanishes unless $W = Y$,
the entries in the above cocomplex 
vanish in all degrees except for
that in which $k\FF_Y$ appears. This degree is precisely
$\codim_X(Y) = \dim(X) - \dim(Y)$, in which case $\Hom_{k\FF}(k\FF_Y,k_Y)
\cong k$.
\end{proof}

\begin{Corollary}
\label{Corollary: quiver of face semigroup algebra}
The quiver $\QQ$ of $k\FF$ is given by the Hasse
diagram of the intersection lattice $\LL$. The cover relations
are oriented by $X \to Y \iff X \lcover Y$.
\end{Corollary}
\begin{proof}
The vertices of $\QQ$ are in one-to-one
correspondence with the isomorphism classes of simple $k\FF$-modules.
These are indexed by the elements of
$\LL$. The number of arrows $X \to Y$ is 
\begin{align*}
\hspace{0.5in}
\dim \Ext^1_{k\FF}(k_X, k_Y) = 
\begin{cases}
1, & \text{if } X \lcover Y, \\
0, & \text{otherwise}.
\end{cases}
\hspace{0.5in} \qedhere
\end{align*}
\end{proof}

\subsection{Quiver Relations}
\label{subsection: quiver relations}

This section defines a $k$-algebra surjection $\varphi: k\QQ \to k\FF$
and identifies a minimal generating set of the kernel. The kernel is
an admissible ideal of the path algebra $k\QQ$, so this generating set
gives the quiver relations.

\subsubsection*{\ref{subsection: quiver relations}A. First Version}
Let $\partial: k\FF \to k\FF$ be the map
$$
\partial(y) = \sum_{x \in k\FF \atop x \gtrdot y} [y:x] x,
$$
where $[y:x]$ is the incidence number defined in equation
(\ref{incidence number identity}).
Define a $k$-algebra morphism $\varphi: k\QQ \longrightarrow k\FF$ by
\begin{gather*}
\varphi(X) = e_X \text{ for } X \in \QQ_0, \\
\varphi(X\to Y) = e_Y \partial(y) e_X, \\
\varphi(X_0 \to X_1 \to \cdots \to X_r) = 
  \varphi(X_{r-1} \to X_r) \cdots \varphi(X_0 \to X_1),
\end{gather*}
where $y$ was chosen in the construction of $e_Y$. (Actually, $y$
can be any element of support $Y$. This follows from the identity
$xx' = x$ iff $\supp(x) \geq \supp(x')$.)
Using Lemma \ref{idempotent lemma} and that
$e_Y = y - \sum_{Z > Y} y e_Z$, it follows that
$e_Y \partial(y) e_X = ([y:x_1] x_1 + [y:x_2]x_2)e_X$
where $x_1$ and $x_2$ are the two faces of support $X$ with common
codimension one face $y$. In particular, this is nonzero.

\begin{Proposition}
\label{quiver relations}
 Let $\varphi: k\QQ \to k\FF$ be the map defined above.
 For each interval $[Z,X]$ of length two in $\LL$, 
 the sum of all paths of length two from $X$ to $Z$
 $$\sum_{Y \in (Z,X)} (X \to Y \to Z)$$ is an element
 of the kernel of $\varphi$. These elements form a minimal
 generating set of relations for the kernel of $\varphi$.
\end{Proposition}

\begin{proof}
If $R$ is a minimal set of relations generating $\ker\varphi$, then
Proposition \ref{minimal relations} gives that the number of elements of
$Z.R.X$ (the number of relations in $R$ starting at $X$ and ending
at $Z$) is $\dim \Ext^2_{k\FF}(k_X,k_Z)$. This is
1 if $[Z,X]$ is an interval of length two and 0
otherwise. Therefore, we need only one relation for each interval of
length two in $\LL$.

Let $z$ be the element of support $Z$ chosen in the construction of $e_Z$.
Then $\sum_{Y \in (Z,X)} \varphi(X \to Y \to Z)$ is a linear combination of
elements of the form $\tilde x e_X$ with $\tilde x$ of support $X$ having $z$
as a face.
If $\tilde x$ has $z$ as a face, then $z$ is of
codimension two in $\tilde x$. Lemma \ref{incidence number lemma} gives that
$\tilde x$ has exactly two codimension one faces $\tilde y$ and $\tilde w$.
Since
\begin{align*} 
 \varphi(\supp(\tilde y) \to Z) & \varphi(X \to \supp(\tilde y)) \\
 & = ([z:\tilde y] \tilde y + [z: y'] y')([y:x_1]x_1 + [y:x_2]x_2) e_X
\end{align*}
and one of $\tilde y x_1$ or $\tilde y x_2$ must be $\tilde x$ ---
suppose $\tilde y x_1 = \tilde x$ --- we see that $\tilde x e_X$ appears
in $\varphi(X \to \supp(\tilde y) \to Z)$ with coefficient 
$[z: \tilde y] [y: x_1]$. The identity (\ref{incidence number identity})
gives this coefficient is $[z: \tilde y] [\tilde y: \tilde x]$.
Similarly, $\tilde x e_X$ appears in $\varphi(X \to \supp(\tilde w) \to
Z)$ with coefficient $[z: \tilde w] [\tilde w: \tilde x]$.
Lemma \ref{incidence number lemma} shows that these two coefficients
sum to zero. Therefore, $\sum_{Y \in (Z,X)} \varphi(X \to Y \to Z) = 0$.
\end{proof}

\begin{Corollary}
The face semigroup algebra $k\FF$ of a hyperplane arrangement depends only on
the intersection lattice $\LL$. 
\end{Corollary}
Note that this implies that arrangements with the same intersection lattice
but nonisomorphic face posets have isomorphic face semigroup algebras.

\subsubsection*{\ref{subsection: quiver relations}B. Second Version}

In this section we note that the idempotents $e_X$ used in the previous
section to define $\varphi$ can be changed slightly without affecting the
kernel of $\varphi$. This will be used in a subsequent paper to construct
idempotents for the descent algebra of a finite Coxeter group
\cite{Saliola2006:ReflectionArrangements-arxiv}.

For each $X \in \LL$ let $L_X$ denote a nonempty set of elements of support
$X$ and let $\lambda_X = |L_X|$. In what follows we will need that the
characteristic of $k$ does not divide $\lambda_X$ for all $X \in \LL$.  Let
$\tilde X$ denote the sum of the elements in $L_X$ divided by $\lambda_X$.
Then $\tilde X$ is an idempotent and the elements $e_X = \tilde X - \sum_{Y
> X} \tilde X e_Y$ form a complete system of primitive orthogonal
idempotents in $k\FF$ (see Remark \ref{Remark: different idempotents}).
Define $\varphi: k\QQ \to k\FF$ using these idempotents: the image of
vertex $X$ is the idempotent $e_X$; the image of an arrow $X \to Y$ is $e_Y
\partial(y) e_X$, where $y$ is any element of support $Y$.

To see that the kernel of $\varphi$ is described by
Proposition \ref{quiver relations}, 
let $(X \to Y \to Z)$ be a path in $\QQ$
and note that $\varphi(X \to Y \to Z)$ can be written as
\begin{align*}
& \frac{1}{\lambda_Z} 
    \sum_{z \in L_Z}
    \Big( [z: y_1^z] y_1^z +  [z: y_2^z] y_2^z \Big) 
    \Big( [y: x_1^y] x_1^y +  [y: x_2^y] x_2^y \Big) e_X,
\end{align*}
where $y_1^z$ and $y_2^z$ are the two faces of support $Y$ with $z$ as a
face and $x_1^y$ and $x_2^y$ are the two faces of support $X$ with $y$ as a
face. (Use Lemma \ref{idempotent lemma}; that $e_X = \tilde X - \sum_{Y >
X} \tilde X e_Y$ for all $X \in \LL$; and Proposition \ref{supp properties}.)



Next we will show that the ceofficient of $y_i^zx_j^y$ in the above
is $\frac{1}{\lambda_Z} [z:y_i^z][y:x_j^y]$. This amounts to showing
that if $y_i^zx_j^y = y_{i'}^{z'}x_{j'}^{y}$, then $z = z'$, $i = i'$ 
and $j = j'$. Well, both $z$ and $z'$ are faces
of $y_i^zx_j^y = y_{i'}^{z'}x_{j'}^{y}$, but no face can have two distinct
faces of the same support. So $z = z'$. Also, $y_i^z$ and $y_{i'}^z$ are
faces of $y_i^zx_j^y = y_{i'}^{z}x_{j'}^{y}$ of the same support, so in
fact $i = i'$.  Since $y_i^zx_j^y = y_{i}^{z}x_{j'}^{y}$, it follows that
$x_j^y$ and $x_{j'}^{y}$ are on the same side of $Y$. But, by definition,
they are on different sides of $Y$. So $j = j'$.

Let $x \in \FF$ have support $X$ and
suppose $x e_X$ is a summand of $\varphi(X \to Y \to Z)$. Then $x = y_i^z
x_j^y$ for some $i,j \in \{1,2\}$, $z \in L_Z$.  Since there are exactly
two faces $w_1$ and $w_2$ in the open interval $\{w \in \FF : z < w < x\}$,
it follows that $y_i^z$ is either $w_1$ or $w_2$. In the former case the
coefficient of $xe_X$ is 
$$\frac{1}{\lambda_Z}[z: w_1][y: x_j^y] = 
\frac{1}{\lambda_Z}[z: w_1][w_1y': w_1x_j^y] = \frac{1}{\lambda_Z}[z:
w_1][w_1: x],$$ using Equation (\ref{incidence number
identity}). Similarly, if $y = w_2$, then the coefficient
is $\frac{1}{\lambda_Z}[z: w_2][w_2: x]$. 
Therefore, the coefficient of $xe_X$ in $\sum_{X \lessdot Y \lessdot Z}
\varphi(X \to Y \to Z)$ is, by Lemma \ref{incidence number lemma},
\begin{align*}
\frac{1}{\lambda_Z}[z: w_1][w_1: x] + \frac{1}{\lambda_Z}[z: w_2][w_2: x]
= 0.
\end{align*}
So $\sum_{X \lessdot Y \lessdot Z} \varphi(X \to Y \to Z) = 0$ since $\{x
e_{X} : \supp(x) = X\}$ is a basis of $k\FF e_X$.

\section{The Ext-algebra of the Face Semigroup Algebra}
 \label{section: the ext-algebra of kF}

\subsection{Koszul Algebras} 
 \label{section: koszul algebras}

Our treatment of Koszul algebras closely follows \cite{BGS1996}.
Let $k$ be a field. A $k$-algebra $A$ is a \emph{graded $k$-algebra} if there
exists a $k$-vector space decomposition $A \cong \bigoplus_{i \geq 0}
A_i$ satisfying $A_iA_j \subset A_{i+j}$. Here
$A_iA_j$ is the set of elements 
$\{ \sum_l a_la'_l \mid a_l \in A_i, a'_l \in A_j \}$. 
The subspace $A_0$ is considered an $A$-module by
identifying it with the $A$-module $A/\bigoplus_{i>0}A_i$.

If $A = \bigoplus_{i\geq0} A_i$ is a graded $k$-algebra, then a \emph{graded
$A$-module} $M$ is an $A$-module with a vector space decomposition $M =
\bigoplus_{i\in\mathbb Z} M_i$ satisfying $A_i M_j \subset M_{i+j}$ for all
$i,j\in\mathbb Z$. A graded $A$-module $M$ is \emph{generated in degree $i$} if
$M_j = 0$ for $j < i$ and $M_j = A_{j-i}M_i$ for all $j \geq i$. 
If $M$ and $N$ are graded
$A$-modules, then an $A$-module morphism $f: M \to N$ has \emph{degree} $p$
if $f(M_i) \subset N_{i+p}$ for all $i$.

A graded $A$-module $M$ has a \emph{linear resolution} if $M$ admits a
projective resolution
\begin{gather*}\begin{CD}
\cdots @>>> P_2 @>d_2>> P_1 @>d_1>> P_0 @>d_0>> M @>>> 0,
\end{CD}\end{gather*}
with $P_i$ a graded $A$-module generated in degree $i$ and $d_i$ a degree 0
morphism form all $i\geq0$. Observe that if $M$ admits a linear resolution,
then $M$ is generated in degree $0$.

\begin{Definition}
A graded $k$-algebra $A = \bigoplus_{i\geq0}A_i$ is a \emph{Koszul algebra}
if $A_0$ is a semisimple $k$-algebra and $A_0$, considered as a graded
$A$-module concentrated in degree $0$, admits a linear resolution.
\end{Definition}

A \emph{quadratic $k$-algebra} is a graded $k$-algebra $A =
\bigoplus_{i\geq0}A_i$ such that $A_0$ is semisimple and $A$ is generated by
$A_1$ over $A_0$ with relations of degree 2. Explicitly, $A =
\bigoplus_{i\geq0}A_i$ is quadratic if $A_0$ is semisimple and $A$ is a
quotient of the free tensor algebra $T_{A_0} A_1 =
\bigoplus_{i\geq0}(A_1)^{\otimes i}$ of the $A_0$-bimodule $A_1$
by an ideal generated by elements of degree 2: $A \cong T_{A_0}A_1/\langle
R\rangle$ with $R \subset A_1 \otimes_{A_0} A_1$.
Here
$(A_1)^{\otimes i}$ denotes the $i$-fold tensor product of $A_1$ over $A_0$.

\begin{Proposition}[\cite{BGS1996}, Corollary 2.3.3]
 Koszul algebras are quadratic.
\end{Proposition}
Not all quadratic algebras are Koszul algebras. Furthermore, it is not known 
for which algebras the notions of quadratic and Koszul coincide.

Let $A = T_{A_0}A_1/\langle R\rangle$ be a quadratic algebra. If $V$ is an
$A_0$-bimodule, let $V^* = \Hom_{A_0}(V,A_0)$. For any subset $W \subset
V$, let $W^\perp = \{ f \in V^* \mid f(W) = 0 \}$.  The algebra $$A^! =
T_{A_0}A_1^*/\langle R^\perp\rangle$$ is the \emph{quadratic dual} of $A$
or the \emph{Koszul dual} of $A$ in the case when $A$ is a Koszul algebra.
(There is an important technicality. In defining the quadratic dual the
identification 
$(V_1^* \otimes \cdots \otimes V_n^*) 
 \cong (V_n \otimes \cdots \otimes V_1)^*$
has been made, where
$(f_1 \otimes \cdots \otimes f_n)(v_n \otimes \cdots \otimes v_1) 
  = f_n(v_n f_{n-1}(v_{n-1} \cdots f_1(v_1)\cdots))$
for all $f_i \in V_i^*$ and $v_i \in V_i$.)

If $A$ is a graded $k$-algebra, then the \emph{$\Ext$-algebra} of $A$ is the
graded $k$-algebra $\Ext(A) = \bigoplus_n \Ext^n(A_0,A_0)$ with
multiplication given by Yoneda composition. 

\begin{Theorem}[\cite{BGS1996}, Theorem 2.10.1 and Theorem 2.10.2]
\label{koszul dual is ext algebra}
 Suppose $A$ is a Koszul algebra. Then the Koszul dual $A^!$ is a Koszul
 algebra isomorphic to the \emph{opposite} of the $\Ext$-algebra 
 $\Ext(A)$ of $A$ and $\Ext(\Ext(A)) \cong A$.
\end{Theorem}

Before proceeding, we record how
the quadratic dual of a quadratic algebra arising as the quotient of the
path algebra of a quiver is constructed from the quiver and relations. 
Note that the path algebra $kQ$ of a quiver $Q$ is the free tensor algebra
of the $k$-vector space $kQ_1$ spanned by the arrows of $Q$ viewed as a
bimodule over the $k$-vector space $kQ_0$ spanned by the vertices of $Q$.
It follows that $A = kQ/\langle R\rangle \cong T_{kQ_0}kQ_1 / \langle R
\rangle$ where $R$ is a set of relations of paths of length two. 
Then the quadratic dual algebra $A^! \cong T_{kQ_0}(kQ_1)^*/ \langle
R^\perp \rangle \cong kQ^{opp}/\langle R^\perp \rangle$ is a quotient of
the path algebra $kQ^{opp}$ of the opposite quiver $Q^{opp}$ of $Q$ and
$R^\perp = \{ s \in kQ^{opp}_2 \mid s^*(r) = 0 \text{ for all } r \in R\}$.
Here $(pq)^*: kQ_2 \to k$ for a path $pq$ of length two in $Q^{opp}$ is the
function that takes the value 1 on $qp \in Q$ and 0 otherwise. 
That is, the quiver of $A^!$ is $Q^{opp}$ and the relations are the
relations 
\emph{orthogonal} to $R$. (This can be derived from the definitions.
See also \cite{GreenMartinez1998}).

%
%
%

\subsection{The Face Semigroup Algebra is a Koszul Algebra}

This section establishes that the face semigroup algebra of a hyperplane
arrangement admits a grading making it a Koszul algebra. This is done by
constructing a linear resolution for the degree 0 component with respect to
the grading inherited from the path length grading on the path algebra of
the quiver.

\begin{Proposition}
 \label{kF is a Koszul algebra}
 $k\FF$ admits a grading making it a Koszul algebra.
\end{Proposition}

\begin{proof}
The $k$-vector spaces 
$$
(k\FF)_i = \bigoplus_{\codim_Y(X)=i} e_X k\FF e_Y.
$$
define a grading on $k\FF$. (This is the grading inherited from the path
length
grading on the path algebra $k\QQ$ of the
quiver $\QQ$ of $k\FF$.) So $k\FF$ is a graded $k$-algebra.
The degree 0 component is
$$
(k\FF)_0 = \bigoplus_{\codim_Y(X)=0} e_X k\FF e_Y
 = \bigoplus_{X \in \LL} e_X k\FF e_X \cong k^{|\LL|},
$$
hence is semisimple. It remains to show that $k^{|\LL|}$ has a linear
resolution. It suffices to show that each simple $k\FF$-module $k_X$ 
has a linear resolution since $k^{|\LL|} \cong \bigoplus_{X\in\LL}
k_X$.

Fix $X \in \LL$ and consider the projective resolution of the simple
$k\FF$-module $k_X$ given by Proposition \ref{projective resolution},
$$
 \cdots \longrightarrow
\left(
 \bigoplus_{\codim_X(Y) = p} k \FF e_Y
\right) \tot{\partial}\longrightarrow
 \cdots \tot{\partial}\longrightarrow
 k \FF e_X \longrightarrow 
 k_X \longrightarrow
 0.
$$
For each $k\FF e_Y$ define $k$-subspaces
$$
(k\FF e_Y)_i 
= \bigoplus_{\codim(W) = i} e_W k\FF e_Y.
$$
By Lemma \ref{idempotent lemma}, if $i < \codim(Y)$, then the degree $i$
component of $k\FF e_Y$ is 0. For $i = \codim(Y)$, $(k\FF)_i = e_Y k\FF e_Y =
\operatorname{span}_k{e_Y}$ (Lemma \ref{idempotent lemma} again). Since $e_Y$
generates $k\FF e_Y$ as a $k\FF$-module, $k\FF e_Y$ is generated in degree
$\codim(Y)$. The boundary operator $\partial$ is a degree 0
morphism: if $e_W w \in e_Wk\FF e_Y$, then $\deg(e_W w) = \codim(W)$ and the
degree of its image $\partial(e_W w) = e_W \partial(w) \in e_W
\partial(k\FF e_Y) \subset \bigoplus_{\codim_X(Y') = p} e_W k\FF e_{Y'}$ is
$\codim(W)$.
\end{proof}

\begin{Remark}
Notice that in creating the surjection $\varphi:k\QQ \to k\FF$ many choices
were taken (in constructing the complete system of primitive orthogonal
idempotents and in putting orientations on the subspaces in $\LL$). These
choices affect the grading inherited by $k\FF$ from $k\QQ$, but the
corresponding graded algebras are isomorphic: two gradings on a $k$-algebra
that both give rise to a Koszul algebra give isomorphic graded
$k$-algebras.  See Corollary $2.5.2$ of \cite{BGS1996}.  
\end{Remark}

\subsection{The $\Ext$-algebra of the Face Semigroup Algebra}

In this section we show that the $\Ext$-algebra of $k\FF$ is the incidence
algebra of the opposite lattice $\LL^*$ of the intersection lattice $\LL$. 

The \emph{incidence algebra} $I(P)$ of a finite poset $P$ is the set of
functions on the subset of $P \times P$ of comparable elements $\{ (y,x) \in
P \times P \mid y \leq x \}$ with multiplication $(fg)(x,y) = \sum_{x \leq z
\leq y} f(x,z)g(z,y)$. The identity element is the Kr\"onecker
$\delta$-function. The incidence algebra $I(P)$ is a split basic algebra and
%
%
the quiver $Q$ of $I(P)$ has $P$ as its set of vertices and exactly one arrow
$x \to y$ if $y \rcovers x$. If $I$ denotes the ideal of $kQ$ generated by
differences of parallel paths, then $I(P) \cong kQ/I$. This isomorphism is
given by mapping a vertex $x$ of $Q$ to the function $y \mapsto \delta(x,y)$,
and an arrow $x \to y$ of $Q$ to the function $(u,v) \mapsto \delta(x,u)
\delta(y,v)$.

\begin{Proposition}
 The Ext-algebra of $k\FF$ is the incidence algebra $I(\LL^*)$ of the 
 opposite lattice of the intersection lattice $\LL$. Equivalently, it is the
 opposite algebra $I(\LL)^{opp}$ of the incidence algebra $I(\LL)$ of $\LL$.
\end{Proposition}
\begin{proof}
Since $k\FF$ is a Koszul algebra (Proposition \ref{kF is a Koszul algebra}),
its Ext-algebra is its Koszul dual algebra (Theorem \ref{koszul dual is ext
algebra}), so we compute the Koszul dual of $k\FF$. 

Let $\QQ$ denote the quiver of $k\FF$.
From Proposition \ref{quiver relations}, $k\FF \cong k\QQ/\langle R \rangle$
is the quotient of the path algebra $k\QQ$ by the ideal generated by the sums
of all parallel paths of length two,
\begin{displaymath}
 R = \left\{ 
 \sum_{Z \in (Y,X)} \left( X \to Z \to Y \right) : X,Y \in \LL
 \right\}.
\end{displaymath}
Then $(k\FF)^! \cong k\QQ^{opp}/\langle R^\perp \rangle$ where
$R^\perp$ is spanned by differences of parallel paths of
length two in $\QQ^{opp}$,
$$
 R^\perp = \left\{ 
 \left( X \to Z \to Y \right) - \left( X \to Z' \to Y \right) : 
 X \rcovers Z,Z' \rcovers Y \in \LL
 \right\}.
$$
(See the
discussion at the end of Section \ref{section: koszul algebras}.)

Let $I(\LL^*)$ denote the incidence algebra of $\LL^*$. 
Then $I(\LL^*) \cong k\QQ^{opp}/I$,
where $I$ is the ideal generated by differences of parallel paths (not
necessarily of length two). Therefore, the proof is complete once it is
shown that $R^\perp$ generates $I$.

If 
$p: X \to X_1 \to \cdots \to X_n \to Y$ and 
$q: X \to Y_1 \to \cdots \to Y_n \to Y$ 
are parallel paths in $\QQ$ such that
there exists an $i$ with $X_j = Y_j$ for all 
$j \neq i$, then $p - q \in I$.
If there exists a sequence of paths $p = p_0, p_1, \ldots, p_j = q$ with
$p_{i-1}$ and $p_{i}$ differing in exactly one place for 
$1 \leq i \leq j$, then 
$p - q = (p_0 - p_1) + \cdots + (p_{j-1} - p_j) \in I$.
Therefore, $I = \langle R^\perp \rangle$ if any path in $\QQ^{opp}$ can be
obtained from any other path that is parallel to it by swapping one vertex
at a time (without breaking the path).  This follows from the
semimodularity of $\LL^*$ and by induction on the length of paths in
$\QQ^{opp}$.  Recall that a finite lattice $L$ is \emph{(upper)
semimodular} if for every $x$ and $y$ in $L$, if $x$ and $y$ cover $x
\wedge y$, then $x \vee y$ covers $x$ and $y$. 

Let 
$X \to X_1 \to \cdots \to X_n \to Y$ and 
$X \to Y_1 \to \cdots \to Y_n \to Y$ be parallel paths in
$\QQ^{opp}$.
Since $X_n$ and $Y_n$ cover $X_n \wedge Y_n = Y$, semimodularity 
of $\LL^*$ 
gives that $X_n \vee Y_n$ covers both $X_n$ and $Y_n$.
Since $X \leq X_n$ and $X \leq Y_n$, it follows that $X \leq
(X_n \vee Y_n)$. So there exists a path from $X$ to $X_n \vee Y_n$.
We are now in the following situation.
\begin{gather*}
\xymatrix@R=0.3em{
                  & X_1 \ar[r] & \cdots \ar[r] & 
X_{n-1} \ar[r] & X_n \ar[rd] & \\
X \ar[ru] \ar[rd] \ar[r] & \cdot \ar[r]    & \cdot \ar[r] & 
(X_n \vee Y_n) \ar[ru] \ar[rd]       &             & Y \\
                  & Y_1 \ar[r] & \cdots \ar[r] & 
Y_{n-1} \ar[r] & Y_n \ar[ru] & \\
}
\end{gather*}
Induction on the length of paths gives that
\begin{gather*}
(Y \to \cdots \to  Y_{n-1} \to Y_n \to Y) 
  - (X \to \cdots \to (X_n \vee Y_n) \to Y_n \to Y), \\
(X \to \cdots \to  X_{n-1} \to X_n \to Y)  
  - (X \to \cdots \to (X_n \vee Y_n) \to X_n \to Y)
\end{gather*}
are in $\langle R^\perp \rangle$. Clearly,
\begin{align*}
(X \to \cdots \to & (X_n \vee Y_n) \to X_n \to Y) \\
 & - (X \to \cdots \to (X_n \vee Y_n) \to Y_n \to Y) \in \langle R^\perp
\rangle.
\end{align*}
Therefore,
\begin{gather*}
(Y \to \cdots \to  Y_{n-1} \to Y_n \to Y) 
  - 
(X \to \cdots \to  X_{n-1} \to X_n \to Y)  
\end{gather*}
is in $\langle R^\perp \rangle$.
Therefore, $I = \langle R^\perp \rangle$ and $(k\FF)^! \cong
k\QQ^{opp}/\langle R^\perp \rangle
= k\QQ^{opp}/I \cong I(\LL^*).$
\end{proof}

\begin{Corollary}
 The $\Ext$-algebra of $I(\LL^*)$ is isomorphic to the face semigroup
 algebra $k\FF$.
\end{Corollary}

\subsection{The Hochschild (Co)Homology of the Face Semigroup Algebra}
 \label{section: hochschild cohomology}

Let $A$ be a $k$-algebra and $M$ an $A$-bimodule. There is a complex
of $A$-bimodules
$$
\begin{CD}
\cdots @>d_{i+1}>> M \otimes_k A^{\otimes i}
 @>d_i>> \cdots 
 @>d_1>> M \otimes_k A
 @>d_0>> M
\end{CD}
$$
with maps $d_i: M \otimes_k A^{\otimes i} \to M \otimes A^{\otimes i-1}$
defined by $d_0(m \otimes a) = am - ma$ for $m \in M$, $a \in A$ and
for $i \geq 1$
\begin{align*}
d_i(m \otimes a_1 \otimes \cdots \otimes a_i) 
= & \ (ma_1 \otimes a_2 \otimes \cdots \otimes a_i) \\
& + \sum_{j=1}^{i-1} (-1)^j (m \otimes a_1 \otimes \cdots \otimes a_ja_{j+1}
 \otimes \cdots \otimes a_i) \\
& + (-1)^i(a_im \otimes a_1 \otimes \cdots \otimes a_{i-1}),
\end{align*}
where $m \in M$ and $a_1, \ldots, a_i \in A$. The \emph{Hochschild homology}
of $A$ with coefficients in $M$ is $\HH_i(A,M) =
\ker(d_i)/\im(d_{i+1})$ for $i \geq 0$. 
Let $\HH_i(A) = \HH_i(A,A)$.

Similarly, there exists a cocomplex of $A$-bimodules
$$
\begin{CD}
M @>d^0>> \Hom_k(A, M)
 @>d^1>> \Hom_k(A\otimes_kA,M)
 @>d^2>> \cdots
\end{CD}
$$
where $d^0: M \to \Hom_k(A,M)$ is the map $d^0(m)(a) = am - ma$
and $d^i$ is the map $d^i: \Hom_k(A^{\otimes i}, M) \to \Hom_k(A^{\otimes
i+1}, M)$ given by
\begin{align*}
(d^if)(a_1 \otimes \cdots \otimes a_{i+1}) 
= & a_1 f(a_2 \otimes \cdots \otimes a_{i+1}) \\
& + \sum_{j=1}^i (-1)^j f(a_1 \otimes \cdots \otimes a_ja_{j+1} \otimes
 \cdots \otimes a_{i+1}) \\
& + (-1)^{i+1} f(a_1 \otimes \cdots \otimes a_i) a_{i+1},
\end{align*}
where $f \in \Hom_k(A^{\otimes i},M)$ and $a_1, \ldots, a_{i+1} \in A$. The
\emph{Hochschild cohomology} of $A$ with coefficients in $M$ is 
$\HH^i(A,M) = \ker(d^i)/\im(d^{i-1})$ for $i \geq 0$.
Denote the Hochschild cohomology of $A$ with coefficients in $A$ by $\HH^i(A)
= \HH^i(A,A)$.

\begin{Proposition}
 The Hochschild homology $\HH_i(k\FF)$ and cohomology 
 $\HH^i(k\FF)$ of $k\FF$
 vanish in positive degrees. In degree zero the homology is
 $\HH_0(k\FF) \cong k^{\#\LL}$ and the cohomology is $\HH^0(k\FF)
 \cong k$.
\end{Proposition}
\begin{proof}
Let $\QQ$ denote the quiver of $k\FF$.  The Hochschild homology of algebras
whose quivers have no oriented cycles is known to be zero in positive
degrees and $k^q$ in degree 0, where $q$ is the number of vertices in the
quiver \cite{Cibils1986}. This establishes the Hochschild homology of
$k\FF$ since $\QQ$ has no oriented cycles. 

Buchweitz (\S3.5 of \cite{Keller2003}) proved that the Hochschild
cohomology algebra of a Koszul algebra is the Hochschild cohomology algebra
of its Koszul dual. Since $k\FF$ is a Koszul algebra with Koszul dual the
incidence algebra $I(\LL^*)$ of the lattice $\LL^*$, there is an
isomorphism
$$
 \HH^*(k\FF) \cong \HH^*(I(\LL^*)) 
 \cong \bigoplus_{i \geq 0} \HH^i(I(\LL^*)).
$$
Gerstenhaber and Schack (\cite{GerstenhaberSchack1983}; also see
\cite[Corollary 1.4]{Cibils1989})
proved that the Hochschild cohomology $\HH^i(I(\LL^*))$ of $I(\LL^*)$ is the
simplicial cohomology of the simplicial complex $\Delta(\LL^*)$ whose
$i$-simplices are the chains of length $i$ in the poset $\LL^*$.
Therefore,
$$
 \HH^i(I(\LL^*)) \cong H^i(\Delta(\LL^*), k).
$$
The latter is zero in positive degrees since $\Delta(\LL^*)$ is a double
cone ($\LL^*$ contains both a top and bottom element) and is $k$
in degree zero since $\Delta(\LL^*)$ is connected. It is easy to check
directly
that $\HH^0(k\FF) \cong k$, completing the proof.
\end{proof}

\section{Connections with Poset Cohomology}
 \label{section: connections with poset cohomology}

\subsection{The Cohomology of a Poset}

Let $P$ denote a finite poset. The \emph{order complex} $\Delta(P)$ of $P$
is the simplicial complex with $i$-simplices the chains of length $i$ in
$P$. Suppose $P$ has both a minimal element $\hat 0$ and a maximal
element $\hat 1$ and let $k$ denote a field. 
The \emph{order cohomology of $P$ with coefficients in $k$} is the
reduced simplicial cohomology with coefficients in $k$ of the order complex
$\Delta(P - \{\hat 0, \hat 1\})$ of $P - \{\hat 0, \hat 1\}$. The order
cohomology of $P$ has the following characterization in terms of the chains
of $P$.

Suppose $P$ contains at least two distinct elements.
For $i \geq 0$, let $C_i(P)$ denote the $k$-vector space
spanned by the $i$-chains of $P - \{ \hat 0, \hat 1 \}$,
$$
C_i(P) = \operatorname{span}_k 
 \Big\{ (x_0 < \cdots < x_i) \ \Big|\ x_j \in P - \{\hat0,\hat1\}\Big\}.
$$
For $i = -1$, let $C_{-1}(P) = k$, the vector
space spanned by the empty chain.
If $P$ consists of one element, then define $C_{-2}(P) = k$ and 
$C_{i}(P) = 0$ otherwise. 

Define coboundary morphisms $\delta_i: C_i(P) \to C_{i+1}(P)$ by
\begin{align*}
 \delta_i(x_0 & < \cdots < x_i) \\ 
  & = \sum_{j=0}^{i+1} (-1)^j \sum_{ x_{j-1} < x < x_j }
    (x_0 < \cdots < x_{j-1} < x < x_j < \cdots < x_i),
\end{align*}
where $x_{-1} = \hat 0$ and $x_{i+1} = \hat 1$. It is straightforward to
check that $\delta^2 = 0$.
The order cohomology of $P$ is $H^i(P) = H^i(P; k) =
\ker(\delta_i)/\operatorname{im}(\delta_{i-1})$.

Notice that if $P$ consists of exactly one element, then 
$H^{-2}(P) = k$ and $H^i(P) = 0$ for $i \neq -2$.
If $P = \{\hat0,\hat1\}$, then 
$H^{-1}(P) = k$ and $H^i(P) = 0$ for $i \neq -1$.

\subsection{A Vector Space Decomposition of the Face Semigroup Algebra}

Suppose the length of the longest chain in the poset $P$ is $d+2$.
Then
$\ker(\delta_{d})$ is spanned by the chains of length $d$ in $P -
\{\hat0,\hat1\}$ and $\operatorname{im}(\delta_{d-1})$ is
spanned by the elements, 
\begin{gather} 
 \label{top cohomology relations}
  \sum_{ x_{j-1} < x < x_j }
    (x_0 \lessdot \cdots \lessdot x_{j-1} \lessdot x \lessdot x_j \lessdot
\cdots \lessdot x_{d-1}),
\end{gather}
one for each chain 
$x_0 \lessdot \cdots \lessdot x_{j-1} < x_j \lessdot \cdots \lessdot
x_{d-1}$ of length $d-1$.

Put $P = \LL$ in the above and identify the cover relations
with the arrows in $\QQ$. Then the top cohomology of $\LL$ corresponds to
the quotient of the span of the maximal paths in $\QQ$ by the
quiver relations. This gives a vector space isomorphism
$e_{\hat0} k\FF e_{\hat1} \cong H^{d-2}(\LL)$, where the length of
the longest chain in $\LL$ is $d$. 
Folkman \cite{Folkman1966} showed that the cohomology of a geometric
lattice is non-vanishing only in the top degree. Since $\LL^*$ is
a geometric lattice and 
$\Delta(\LL^*) = \Delta(\LL)$,
the cohomology of $\LL$
is non-vanishing only in the top degree. Therefore,
$e_{\hat0} k\FF e_{\hat1} \cong H^*(\LL)$.
Since every interval of a geometric lattice is also a geometric
lattice, the result holds for every interval of $\LL$. 
That is, $e_X k\FF e_Y \cong H^*([X,Y])$.
\begin{Proposition}
 \label{proposition: vector space decomposition into cohomology groups}
$k\FF$ has a $k$-vector space decomposition
in terms of the order cohomology of the intervals of $\LL$,
$$k\FF \cong \bigoplus_{X,Y \in \LL} H^*([X,Y]).$$
\end{Proposition}

\subsection{Another Cohomology Construction on Posets}
 \label{subsection: new cohomology construction}

In light of the above decomposition,
the direct sum $\bigoplus_{X,Y \in \LL} H^*([X,Y])$
inherits a $k$-algebra structure from $k\FF$.  This section shows that
the algebraic structure can be obtained via the cup product of a cohomology
algebra on the intersection lattice. This cohomology construction appears to be
new.

Let $P$ be a finite poset and let $k$ denote a field.  Let $D_i(P)$
denote the $k$-vector space of $i$-chains in $P$,
$$
D_i(P) = \Big\{ (x_0 < \cdots < x_i) \ \Big|\ x_j \in P \Big\}.
$$
Define coboundary morphisms $d_i: D_i(P) \to D_{i+1}(P)$ by
\begin{align*}
 d_i(x_0 & < \cdots < x_i) \\
  & = \sum_{j=1}^{i} (-1)^j \sum_{ x_{j-1} < x < x_j }
    (x_0 < \cdots < x_{j-1} < x < x_j < \cdots < x_i).
\end{align*}
Then $d^2 = 0$. The cohomology groups of the cocomplex 
$(D_\bullet, d)$ will be denoted by
$\FH^i(P) 
 = \FH^i(P; k) = \ker(d_i) / \operatorname{im}(d_{i-1})$.

The differences between $\FH^i(P)$ and $H^i(P)$ are small, but important.
The former is defined for any poset $P$, not just a poset
with $\hat0$ and $\hat1$. The vector space $D_i(P)$ is spanned by
\emph{all} the chains in $P$, not just those avoiding $\hat 0$ and $\hat 1$.
The summation in the
coboundary morphism $d_i: D_i(P) \to D_{i+1}(P)$ runs from $j=1$ to $j=i$,
whereas the summation runs from $j=0$ to $j=i+1$ in the coboundary
morphism
$\delta_i: C_i(P) \to C_{i+1}(P)$. However, there is a strong relationship
between $\FH(P)$ and $H(P)$.

\begin{Proposition}
 \label{proposition: full cohomology relation with order cohomology}
 Let $P$ be a finite poset. Then for all $i \geq 0$,
 $$
  \FH^i(P) 
  \cong \bigoplus_{x, y \in P} \OH^{i-2}([x,y]).
 $$
\end{Proposition}
\begin{proof}
$D_i(P)$ decomposes into subspaces spanned by the $i$-chains of $P$ beginning
at $x$ and terminating at $y$: $(x < x_1 < \cdots < x_{i-1} < y)$. The
differential $d_i$ respects this decomposition and the subspaces are
isomorphic to $C_{i-2}([x,y])$ (drop the $x$ and $y$ of each chain).
This isomorphism commutes with the coboundary operators, establishing the
proposition. 
\end{proof}

The benefit of working with $\FH^*(P)$ is that the 
\emph{simplicial cup product} 
(see \cite[\S49]{Munkres1984}) on the
simplices of the order complex $\Delta(P)$ of $P$ descends to a product on
the cohomology. 

Define a product $\smile: D_p(P) \times D_q(P) \to D_{p+q}(P)$ by
\begin{align}
 \label{cup product}
(x_0 < \cdots < x_p) & \smile (y_0 < \cdots < y_q) \notag\\
 & = \begin{cases}
(x_0 < \cdots < x_p = y_0 < \cdots < y_q), & x_p = y_0, \\
0, & x_p \neq y_0.
\end{cases}
\end{align}
\begin{Lemma}
For $c \in D_p(P)$ and $d \in D_q(P)$, 
\begin{gather*}
\delta_{p+q}(c \!\smile\! d) 
 \hspace{0.5em} = \hspace{0.5em} 
\delta_p(c)\!\smile\! d  
 \hspace{0.5em} + \hspace{0.5em}    
 (-1)^p c \!\smile\! \delta_q(d) .
\end{gather*}
\end{Lemma}
\begin{proof}
Let $c = (x_0 < \cdots < x_p)$ and $d = (x_p < \cdots < x_{p+q})$.
Then
\begin{align*}
& \delta_p(c) \!\smile\! d \ + \  (-1)^p c \!\smile\! \delta_q(d) \\
&= \sum_{j=1}^p (-1)^j \sum_{x_{j-1} < x < x_j}
    ( x_0 < \cdots < x_{j-1} < x < x_j < \cdots < x_{p} ) \!\smile\! d \\
&\hspace{10pt} + (-1)^p c \!\smile\! \hspace{-5pt} 
    \sum_{j=p+1}^{p+q} \hspace{-5pt} (-1)^{j-p}\hspace{-15pt}
    \sum_{x_{j-1} < x < x_j} \hspace{-15pt}
    ( x_p < \cdots < x_{j-1} < x < x_j < \cdots < x_{p+q} ) \\
&= \sum_{j=1}^{p} (-1)^j \sum_{x_{j-1} < x < x_j}
    ( x_0 < \cdots < x_{j-1} < x < x_j < \cdots < x_{p+q}) \\
&\hspace{10pt} + \sum_{j=p+1}^{p+q} (-1)^{j} \sum_{x_{j-1} < x < x_j}
    ( x_0 < \cdots < x_{j-1} < x < x_j < \cdots < x_{p+q} ) \\
&= \sum_{j=1}^{p+q} (-1)^j \sum_{x_{j-1} < x < x_j}
    ( x_0 < \cdots < x_{j-1} < x < x_j < \cdots < x_{p+q} ) \\
&= \delta_{p+q}(x_0 < \cdots < x_{p+q}) \\
&= \delta_{p+q}(c \!\smile\! d).
\qedhere
\end{align*}
\end{proof}

\begin{Corollary}
The product $D_p(P) \times D_q(P) \tot{\smile}\longrightarrow D_{p+q}(P)$ 
induces a well-defined product 
$ \FH^p(P) \times \FH^q(P) \tot\smile\longrightarrow \FH^{p+q}(P)$
giving $\FH^*(P) = \bigoplus_{i} \FH^i(P)$ a $k$-algebra structure.
\end{Corollary}

\subsection{The Face Semigroup Algebra as a Cohomology Algebra}

Combining Propositions \ref{proposition: vector space decomposition into
cohomology groups} and \ref{proposition: full cohomology relation
with order cohomology} 
gives the vector space isomorphism 
\begin{gather*}
 \phi: \FH^*(\LL) 
 \tot{\cong}\to \bigoplus_{X,Y\in\LL} H^*([X,Y])
 \tot{\cong}\to k\QQ/I
 \tot{\varphi}\to k\FF.
\end{gather*}
Recall that Proposition \ref{proposition: vector space decomposition into
cohomology groups} identifies $\bigoplus_{X,Y} H^*([X,Y])$ with 
$k\FF$ via the quiver $\QQ$ with relations of $k\FF$.
The isomorphism identifies an unrefinable chain
in $\LL$ with the corresponding path in $\QQ$
$$
\big(X_0 \lessdot X_1 \lessdot \cdots \lessdot X_{j-1} \lessdot X_j\big)
\longmapsto
\big(X_j \to X_{j-1} \to \cdots \to X_1 \to X_0\big)
$$
and maps the relations in $\FH^*(\LL)$ to the quiver relations.
Under this isomorphism the multiplication in $\FH^*(\LL)$ maps
to the multiplication in $k\QQ/I$ (composition of chains
in $\LL$ maps to composition of paths in $\QQ$).
Therefore, $\phi$ is a $k$-algebra isomorphism.

\begin{Proposition}
Let $k\FF$ be the face semigroup algebra of a hyperplane arrangement
with intersection lattice $\LL$. Then $k\FF \cong \FH^*(\LL)$.
\end{Proposition}

\subsection{Connection with the Whitney Cohomology of the Intersection
Lattice}
We finish this section by identifying the Whitney cohomology of $\LL$
in $k\FF$. 
(See  \cite{Baclawski1975} and more recently \cite{Wachs1999}.)
The \emph{Whitney cohomology} of a poset $P$ with
$\hat0$ is the direct sum
$
\WH^*(P) = \bigoplus_{X \in P} H^*([\hat0, X]).
$
Since the Whitney homology of $\LL^*$ is isomorphic to the Orlik-Solomon
algebra of $\LL^*$
(\cite[\S7.10]{Bjorner1992}), the following result also explains how the
dual of the Orlik-Solomon algebra embeds in the face semigroup algebra.
\begin{Corollary}
The Whitney cohomology of $\LL^*$ is isomorphic to the ideal of chambers
in $k\FF$. It is a projective indecomposable $k\FF$-module.
\end{Corollary}
\begin{proof}
Since $H^*([X,Y]) \cong e_X k\FF e_Y$ for all $X, Y \in \LL$
(see the discussion preceeding Proposition 
\ref{proposition: vector space decomposition into cohomology groups}),
the Whitney cohomology of $\LL^*$ is
\begin{align*}
 \hspace{2em}
 \WH^*(\LL^*) 
 \cong \bigoplus_{X \in \LL} H^*([X,\hat1]) 
 \cong
   \bigoplus_{X \in \LL} e_X k\FF e_{\hat1}
 \cong k\FF e_{\hat1} \cong k\FF_{\hat1}.
 \hspace{2em}\qedhere
\end{align*}
\end{proof}

\section{Future Directions}

These results extend to the semigroup algebra of the semigroup of covectors
of an oriented matroid (see \cite[\S4.1]{OrientedMatroids1993} for the
definition of this semigroup).  This is essential due to two observations.
The first observation is that the exact sequence used to construct the
projective resolutions of the simple modules (Section \ref{section:
projective resolutions of simple modules}) can be extended to the semigroup
algebra of an oriented matroid \cite[\S6]{BrownDiaconis1998}.  The second
observation is that the construction of the complete set of primitive
orthogonal idempotents in Section \ref{primitive idempotents} holds for a
larger class of semigroups call \emph{left regular bands} (see
\cite{Saliola2006:LRBs-arxiv}).

By restricting attention to the reflection arrangement of a finite Coxeter
group, the theory developed here yields results about the descent algebra of
the Coxeter group. In a subsequent paper
\cite{Saliola2006:ReflectionArrangements-arxiv} we will study the quiver
and module structure of the descent algebra using this approach.

The cohomology construction introduced in Section \ref{subsection: new
cohomology construction} appears to be new.  This construction is
interesting, especially because the resulting cohomology algebra appears
naturally as the face semigroup algebra of a hyperplane arrangement, and
deserves to be studied further. The natural starting point would be to
mimic the theory of the order cohomology of a poset. We mention one
possibility: if $G$ is a group acting on a poset $P$, then the $G$-action
on $P$ induces a $G$-module structure on $\FH^*(P)$, and the resulting
$G$-module structure can be studied. This has been extensively studied for
order homology and cohomology and is quite interesting (\cite{Wachs1999},
for example).

For certain classes of posets $\FH^*(P)$ has nice algebraic structure.
For example, if $P$ is a Cohen-Macaulay poset, then its incidence algebra
$I(P)$ is a Koszul algebra \cite{Polo1995}, \cite{Woodcock1998}. Hence,
$\FH^*(P)$ is the Koszul dual algebra of $I(P)$. This describes the Koszul
dual algebra of $I(P)$ in terms of the order cohomology of $P$.

The construction also provides an extension of a result \cite{Hozo1996}
describing a part of the Lie algebra (co)homology of a certain subalgebra
$N(P)$ of the incidence algebra of $P$ in terms of the order (co)homology
of $P$.  Hozo showed that if $P$ contains $\hat 0$ and $\hat 1$, then the
Lie algebra (co)homology of $N(P)$ contains the order (co)homology of $P$.
His proof extends to show that for any poset $P$ (not necessarily
containing $\hat0$ and $\hat1$), the Lie algebra (co)homology of $N(P)$
contains the (co)homology $\FH^*(P)$. This is a further step towards
describing the complete Lie algebra (co)homology of $N(P)$ in terms of the
combinatorics of the poset $P$.

\bibliographystyle{hapalike}
\bibliography{references} 

\end{document}